\documentclass[12pt]{amsart}

\usepackage{amssymb}

\usepackage{enumitem}

\usepackage{graphicx}

\makeatletter
\@namedef{subjclassname@2020}{%
  \textup{2020} Mathematics Subject Classification}
\makeatother

\usepackage[T1]{fontenc}

\newtheorem{theorem}{Theorem}[section]
\newtheorem{corollary}[theorem]{Corollary}
\newtheorem{lemma}[theorem]{Lemma}

\theoremstyle{definition}
\newtheorem{definition}[theorem]{Definition}

\numberwithin{equation}{section}

\newcommand{\afp}{\mathrm{AFP}}
\newcommand{\afs}{\mathrm{AFS}}
\newcommand{\aht}{\mathrm{AHT}}
\DeclareMathOperator{\cf}{cf}
\DeclareMathOperator{\supp}{supp}
\newcommand{\id}{\mathrm{id}}

\begin{document}


\baselineskip=17pt


\title[Uncountable Adjacent Hindman's Theorem]{The Adjacent Hindman's theorem for uncountable groups}

\author[L. Carlucci]{Lorenzo Carlucci}
\address{Dipartimento di Matematica ``G. Castelnuovo"\\ Sapienza Università di Roma\\
00185 Roma, Italy}
\email{lorenzo.carlucci@uniroma1.it}

\author[D. Fern\'andez]{David Fern\'andez-Bret\'on}
\address{Escuela Superior de F\'{\i}sica y Matem\'aticas\\
Instituto Polit\'ecnico Nacional\\
07738, Ciudad de M\'exico, Mexico. }
\urladdr{https://dfernandezb.web.app}
\email{dfernandezb@ipn.mx}



\date{}

\begin{abstract}
Recent results of Hindman, Leader and Strauss and of the second author and Rinot
showed that some natural analogs of Hindman's Theorem fail for all uncountable cardinals. 
Results in the positive direction were obtained by Komj\'ath, the first author, and the second author and Lee, who 
showed that there are arbitrarily large Abelian groups satisfying some Hindman-type property.
Inspired by an analogous result studied by the first author in the countable setting, we prove a new variant of Hindman's Theorem for uncountable cardinals, called the
Adjacent Hindman's Theorem:
For every $\kappa$ there is a $\lambda$ such that, whenever a group $G$ of cardinality $\lambda$
is coloured with $\kappa$ colours, there exists a $\lambda$-sized injective sequence of elements of $G$
with all finite products of adjacent terms of the sequence of the same colour. We obtain bounds on $\lambda$ as a function of $\kappa$, and prove that such bounds are optimal. This is the first example of a Hindman-type result for uncountable cardinals that we can prove also in the non-Abelian setting and, furthermore, it is the first such example where monochromatic
products (or sums) of {\it unbounded length} are guaranteed.
\end{abstract}

\maketitle

\section{Introduction}

Hindman's famous Finite Sums Theorem \cite{Hin:74} states that for any finite colouring of the
positive integers there is an infinite set $H$ of positive integers such that all non-empty finite sums of distinct elements
of $H$ are coloured the same. Hindman's Theorem can be generalized, via the use of ultrafilters, utilizing the Galvin--Glazer argument, to arbitrary right-cancellative infinite semigroups (for further generalizations, see \cite{Gol-Tsa:13}).

The question whether analogues of Hindman's Theorem hold for uncountable cardinals 
received substantial attention in recent years (\cite{Kom:15, SW:15, HLS:16, Bre:16, Bre-Rin:16, Kom:16, Car:19:order, fernandez-lee})\footnote{For a related but different question, going back to Erd\H{o}s~\cite{Erd:73}, see \cite{Ele-Haj-Kom:91}.}. 
The variants of interest range from the full analogue of Hindman's Theorem to restrictions thereof where one varies either the number of colours
or the set of finite sums that are required to be monochromatic.

A series of negative results showed that, in some sense, Hindman's Theorem is only a countable phenomenon. Hindman, Leader and Strauss~\cite{HLS:16} exhibited a colouring of $\mathbb{R}$ in $2$ colours
such that no continuum-sized set has monochromatic $2$-terms sums.
Komj\'ath~\cite{Kom:15} and, independently, Soukup and Weiss \cite{SW:15} improved this by ruling out solutions of uncountable cardinality. 
The second author~\cite{Bre:16} showed that for any uncountable Abelian group $G$ there exists a colouring of $G$ in $2$ colours
such that no uncountable subset of $G$ has all its finite sums of the same colour\footnote{We later found out that this theorem had been proved previously not once, but on two occasions, by Malykhin~\cite[p. 70, within the proof of Theorem 2]{malykhin-extremally-disconnected} (for Boolean groups) and Protasov~\cite[Theorem 2]{protasov-ultrafilters-close-to-ramsey} (for all Abelian groups), although it seems not many people in the mathematical community are aware of these earlier results.}. Further negative results on colourings of uncountable Abelian groups were proved by the second author and Rinot in a 
follow-up paper~\cite{Bre-Rin:16}, where, e.g., the following theorem is established: 
Every uncountable Abelian group can be coloured with countably many colours in such a way that every set of finite sums generated by uncountably many elements must be panchromatic, i.e., must contain elements of all possible colours.

On the positive side, some natural variants of Hindman's Theorem have been proved to hold for uncountable cardinals. 
Komj\'ath~\cite{Kom:16} gave the first such example, showing that {\em there exist} arbitrarily large Abelian groups satisfying 
{\em some form} of Hindman-type theorem. 
In particular: For every finite $n$ and infinite cardinal $\kappa$ there is an Abelian group $G$
such that for every colouring of $G$ with $\kappa$ colours there are $n$ elements such that
the same colour class contains all subsums of them. Moreover, for every cardinal $\kappa$, for every $n>1$,
there exists a sufficiently large Abelian group $G$ such that for every colouring of $G$
with $\kappa$ colours there are distinct elements $\{a_{i,\alpha} : 1 \leq  i \leq  n, \alpha < \kappa\}$ 
such that all sums of the form $a_{i_1,\alpha_1} +\dots + a_{i_r,\alpha_r}$ (for arbitrary distinct values of $i_1, i_2,\dots, i_r\in [1,n]$ and $\alpha_1,\ldots,\alpha_r<\kappa$) 
are distinct and in the same colour class. Note that Komj\'ath's results are existential and are witnessed by a Boolean group of suitably large
size. The second author and Lee~\cite{fernandez-lee} generalized this result to {\em all} sufficiently large Abelian groups in the case $n=2$, and showed that no such generalization is possible for $n\geq 3$.
The first author \cite{Car:19:order} isolated a family of natural restrictions of Hindman's Theorem all of which are satisfied
by every sufficiently large uncountable Abelian group. These restrictions relax the monochromaticity condition to those sums whose length belongs to {\em some
structured set} of positive integers drawn from an uncountable solution set.  
A typical member of this family is the following: 
For any infinite cardinal $\kappa$ and positive integers $c$ and $d$, there exists $\lambda$ such that
{\em for any} Abelian group $G$ of size $\lambda$, for every colouring of $G$ with $c$ colours there is 
a $\lambda$-sized $X\subseteq G$ and positive integers $a,b$ such that all sums of $n$-many distinct elements of $X$ with 
$n \in \{a, a+b, a+2b, \dots, a+db\}$ have the same colour.

In this note we introduce a new natural restriction of Hindman's Theorem for uncountable groups. The idea comes from \cite{Car:18:weak}, where the following so-called Adjacent Hindman's Theorem is introduced in the countable setting: Whenever the positive integers are finitely coloured, there exists an infinite set such that all sums of adjacent elements of that set have the same colour. By sums of adjacent elements of a set $X=\{ x_0, x_1, x_2, \dots\}$ enumerated in increasing order we mean all sums of the form $x_i + x_{i+1} + \dots + x_{i+\ell}$ for $i, \ell$ non-negative integers. This restriction of Hindman's Theorem is interesting from the point of view of Reverse Mathematics, because it has a significantly weaker strength than the full Hindman's Theorem, and in fact it follows from Ramsey's Theorem for pairs over the axiom system $\mathsf{RCA}_0$ formalizing computable mathematics\footnote{Further recent, as of yet unpublished, work of the second author shows that this restriction of Hindman's Theorem is in fact equivalent, over $\mathsf{RCA}_0$, to a very natural weakening of Ramsey's Theorem where one only considers colourings that are translation-invariant.} (the precise strength of Hindman's Theorem is an important and longstanding open problem in Reverse Mathematics, see \cite{Car:2021} for details). We formulate and prove a suitable analogue for uncountable groups and give optimal bounds on the size of a set satisfying the principle as a function of the number of colours. 

This is the first positive example of an uncountable Hindman-type theorem that guarantees monochromatic sums (or products) of {\em unbounded length}. Furthermore, this is also the first case of an uncountable Hindman-type theorem that holds not just of Abelian groups, but of all groups.

\section{Adjacent Finite Products and the uncountable Adjacent Hindman's Theorem}

We begin by introducing the terminology needed to state our uncountable Adjacent Hindman's Theorem, and prove a preliminary result to rule out trivial cases. We use multiplicative notation for groups since our groups are not necessarily Abelian. We first generalize the idea of adjacent finite products from countable to uncountable. 

\begin{definition}
If $G$ is a group, $\mu$ is an ordinal, and $(g_\alpha\big|\alpha<\mu)$ is a sequence of elements of $G$, we define its set of {\bf adjacent finite products} to be
\begin{equation*}
\afp(g_\alpha\big|\alpha<\mu)=\{g_\alpha g_{\alpha+1}\cdots g_{\alpha+n}\big|\alpha<\mu\text{ and }n<\omega\}.
\end{equation*}
(On occasion, if the group is Abelian, we may write $\afs$ instead and use additive notation.)
\end{definition}

The products considered are only of finitely many adjacent indices. As a result of that, if we decompose $\mu=\omega\cdot\nu+m$ by the division algorithm, it is always the case that
\begin{equation*}
\afp(g_\alpha\big|\alpha<\mu)=\left(\bigcup_{\xi<\nu}\afp(g_{\omega\cdot\xi+n}\big|n<\omega)\right)\cup\afp(g_{\omega\cdot\nu+n}\big|n<m),
\end{equation*}
and if $\mu$ is limit (i.e., $m=0$ in the above decomposition) then the last term of the above union is empty. Therefore, the set of adjacent finite products of a transfinitely long sequence of elements is simply the union of adjacent finite products sets of many $\omega$-sequences. Thus, we could conceivably take an $\omega$-sequence and then just ``repeat it many times'' to get arbitrarily long sequences with the same $\afp$-sets. To avoid these trivialities, we will focus on injective sequences.

\begin{definition}\hfill
\begin{enumerate}
\item If $G$ is a group, $\mu$ is an ordinal, and $\kappa$ is a non-zero cardinal, we say that $G$ {\em satisfies the Adjacent Hindman's Theorem for $\mu$ and $\kappa$} if for every colouring $c:G\longrightarrow\kappa$ there exists an injective sequence $(g_\alpha\big|\alpha<\mu)$ of elements of $G$ such that $\afp(g_\alpha\big|\alpha<\mu)$ is monochromatic for $c$. We denote this by $G\rightarrow(\mu)_\kappa^\afp$ or $\aht(G, \mu, \kappa)$.

\item The statement $\lambda\rightarrow(\mu)_\kappa^\afp$ means that for every group $G$ with $|G|=\lambda$, the statement $G\rightarrow(\mu)_\kappa^\afp$ holds. In this case we also say that {\em $\lambda$ satisfies the Adjacent Hindman's Theorem for $\kappa$-colourings} and write $\aht(\lambda,\mu,\kappa)$. 
\end{enumerate}
\end{definition}

Clearly, if $\mu'\leq\mu$ and $\kappa'\leq \kappa$, then $\aht(\lambda,\mu,\kappa)$ implies $\aht(\lambda,\mu',\kappa')$; furthermore, $\aht(\lambda,\mu,\kappa)$ is false (for any $\kappa$) as soon as $\mu>\lambda$. As we will see throughout the paper, for a fixed $\kappa$, most sufficiently large cardinals $\lambda$ make the statement $\aht(\lambda,\lambda,\kappa)$ true (the observation above about decomposing sets of adjacent finite products will be pivotal for that). On the other hand, the principle $\aht(\lambda,2,\kappa)$ is false as soon as $\kappa\geq\lambda$. Indeed, if $|G|=\lambda\leq\kappa$, then one verifies that $G\nrightarrow(2)_\kappa^\afp$ by simply colouring each element of $G$ with a different colour.


Let $G$ be a group. 
A sequence $(g_\alpha\big|\alpha<|G|)$ of elements of $G$ satisfying $g_\alpha\notin\langle g_\xi\big|\xi<\alpha\rangle$ for all $\alpha < |G|$, where $\langle g_\xi\big|\xi<\alpha\rangle$ denotes the subgroup generated by the $g_\xi$'s with $\xi < \alpha$, is called an {\em independent sequence} in $G$. If $G$ is uncountable, then such a sequence can be obtained, by an easy transfinite induction, as in~\cite[Lemma 2]{fernandez-lee}, which is stated for Abelian groups but works in general. 

The following will be extremely useful in our proofs. 

\begin{lemma}\label{lem1}
Let $G$ be a group, and let $(g_\alpha \big| \alpha <|G|)$ be an independent sequence in $G$. If $\alpha, \beta, \gamma, \delta < |G|$ are four ordinals such that either
\begin{enumerate}
\item the ordinals $\alpha,\beta,\gamma,\delta$ are all pairwise distinct, or
\item $\alpha<\beta=\gamma<\delta$,
\end{enumerate}
then $g_\alpha g_\beta^{-1}\neq g_\gamma g_\delta^{-1}$.
\end{lemma}

\begin{proof}
Suppose, aiming for a contradiction, that 
\begin{equation}\label{eqn1}
g_\alpha g_\beta^{-1}=g_\gamma g_\delta^{-1},
\end{equation}
and consider the ordinal $\varepsilon=\max\{\alpha,\beta,\gamma,\delta\}$. If $\alpha,\beta,\gamma,\delta$ are pairwise distinct, then we may assume without loss of generality that $\varepsilon$ is either $\gamma$ or $\delta$; in the other case, it is plain that $\varepsilon=\delta\notin\{\alpha,\beta,\gamma\}$. In case $\varepsilon=\gamma$ and $\alpha,\beta,\gamma,\delta$ are pairwise distinct, by~(\ref{eqn1}) we obtain $g_\gamma=g_\alpha g_\beta^{-1} g_\delta\in\langle g_\xi\big|\xi<\gamma\rangle$, a contradiction. If, on the other hand, we have $\varepsilon=\delta\notin\{\alpha,\beta,\gamma\}$, then~(\ref{eqn1}) yields $g_\delta=g_\beta g_\alpha^{-1} g_\gamma\in\langle g_\xi\big|\xi<\gamma\rangle$, again a contradiction.
\end{proof}

\section{Cardinals satisfying the Adjacent Hindman's Theorem}

We start by proving that the Adjacent Hindman's Theorem $\aht(\lambda,\lambda,\kappa)$ for $\kappa$-colourings holds for ``most'' sufficiently large (relative to $\kappa$) cardinals $\lambda$.
 
Our first result deals with the case where the number of colours is smaller than the cofinality of the size of the group. In this case, if the size of the group satisfies Erd\H{o}s and Rado's generalization of Ramsey's Theorem ensuring a countable monochromatic set then it also satisfies the Adjacent Hindman's Theorem. We note that also in the countable realm, Ramsey's Theorem for pairs implies the (countable) Adjacent Hindman's Theorem (see \cite[Proposition 1]{Car:18:weak}).

Let us fix some terminology and notation and recall some basic facts. 
If $\kappa$ is a cardinal we denote by $\kappa^+$ its successor. 
Let $\beth_0(\kappa)=\kappa$ and $\beth_{n+1}(\kappa)=2^{\beth_n(\kappa)}$.
The Erd\H{o}s--Rado Theorem (see~\cite{EHMR:84}) is the following statement generalizing 
Ramsey's Theorem to the uncountable setting: 
For every infinite cardinal $\kappa$, for every
integer $n$, the following holds: If the $(n+1)$-tuples of a set $X$ of cardinality $\beth_n(\kappa)^+$ are coloured with $\kappa$ colours then there exists a subset $H$ of $X$ such that $H$ has cardinality 
$\kappa^+$ and such that all $(n+1)$-tuples from $H$ have the same colour. 
We use $\lambda \rightarrow(\mu)^n_\kappa$ to denote the fact that, for every set $X$ of cardinality $\lambda$, 
for every colouring of $[X]^n$ with $\kappa$ colours, there exists a monochromatic $H \subseteq X$ of cardinality
$\mu$. In this notation the Erd\H{o}s--Rado Theorem is written 
$\beth_n(\kappa)^+ \to (\kappa^+)^{n+1}_\kappa$ and Ramsey's Theorem is written $\omega \to (\omega)^2_k$
(for $k < \omega$).

\begin{lemma}\label{lem:ub}
Let $\kappa$ be a cardinal and let $\lambda$ be an uncountable cardinal such that $\lambda \rightarrow (\omega)^2_\kappa$ and $\cf(\lambda) > \kappa$. Then $\lambda \rightarrow (\lambda)^{\afp}_\kappa$.
\end{lemma}

\begin{proof}
Let $G$ be a group of cardinality $\lambda$, and let $(g_\alpha \big | \alpha < \lambda)$ be an independent sequence in $G$. Fix a partition $C=\{ C_\alpha \big | \alpha < \lambda\}$
of $\lambda$ into $\lambda$ cells, each of which has cardinality $\lambda$.

Let $c: G\longrightarrow \kappa$ be a colouring. For each $\alpha < \lambda$ let $d_\alpha : [C_\alpha]^2 \longrightarrow \kappa$ be defined by
$$ d_\alpha (\{\beta, \gamma\}) = c(g_\beta g^{-1}_\gamma)$$
whenever $\beta,\gamma \in C_\alpha$ with $\beta < \gamma$. We let $d: [\lambda]^2 \longrightarrow \kappa$ be the union of the $d_\alpha$'s.

Since $\lambda \rightarrow (\omega)^2_\kappa$ by hypothesis, for each $\alpha < \lambda$ there exists an (increasing) sequence $(\xi^\alpha_n \big | n < \omega)$ of elements of $C_\alpha$ and a colour $i_\alpha < \kappa$ such that
$[\{ \xi^\alpha_n \big | n < \omega\}]^2$ is monochromatic for $d_\alpha$ with colour $i_\alpha$. 
Note that for different $\alpha$'s, the sets $\{ \xi^\alpha_n \big | n < \omega\}$ are disjoint since they belong to different cells of the partition $C$. 

Since $\cf(\lambda) > \kappa$, there exists a $Y\subseteq \lambda$ and an $i < \kappa$ such that $|Y| = \lambda$ and for all $\alpha \in Y$ we have $i_\alpha = i$.

Since $\lambda$ is a cardinal number, we have $\lambda = \omega \cdot \lambda$ (ordinal multiplication). 
Pick an injection $\psi : \lambda \to Y$ and define, for each $\alpha < \lambda$ and $n < \omega$,
$$ x_{\omega \cdot \alpha + n } = g_{\xi^{\psi(\alpha)}_n} g^{-1}_{\xi^{\psi(\alpha)}_{n+1}}.$$

We claim that the sequence $(x_\alpha \big | \alpha < \lambda)$ is injective. So, suppose that $x_{\omega \cdot \alpha + n } = x_{\omega \cdot \beta + m }$, and we aim to show that $\alpha = \beta$ and $n = m$.  

Our assumption is that $ g_{\xi^{\psi(\alpha)}_n} g^{-1}_{\xi^{\psi(\alpha)}_{n+1}} = g_{\xi^{\psi(\beta)}_m} g^{-1}_{\xi^{\psi(\beta)}_{m+1}}$. Note that, if $\alpha\neq\beta$, then the factors on the left-hand side belong to $C_{\psi(\alpha)}$ and the factors on the right-hand side belong to $C_{\psi(\beta)}$, with $C_{\psi(\alpha)}\cap C_{\psi(\beta)}=\varnothing$, so this case is ruled out by part (1) of Lemma~\ref{lem1}. Hence we may assume that $\alpha=\beta$ and, without loss of generality, suppose also that $n\leq m$. Then we have $g_{\xi^{\psi(\alpha)}_n} g^{-1}_{\xi^{\psi(\alpha)}_{n+1}}=g_{\xi^{\psi(\alpha)}_m} g^{-1}_{\xi^{\psi(\alpha)}_{m+1}}$. Since the sequence $g_{\xi^{\psi(\alpha)}_n}$ is strictly increasing in $n$, the only possibility that is not ruled out by part (1) of Lemma~\ref{lem1} is that $m\in\{n,n+1\}$. Part (2) of Lemma~\ref{lem1} rules out the case $m=n+1$, so the conclusion is that we must have $n=m$, and we are done.

Finally, it is easy to show that $\afp(x_\alpha \big | \alpha < \lambda)$ is monochromatic for $c$ with colour $i$ by 
construction. Notice that, for every $n,m<\omega$, we have
\begin{eqnarray*}
c(x_{\omega\cdot\alpha+n} x_{\omega\cdot\alpha+n+1} \cdots x_{\omega\cdot\alpha+n+m}) & 
= & c\left((g_{\xi^{\psi(\alpha)}_n} g^{-1}_{\xi^{\psi(\alpha)}_{n+1}})(g_{\xi^{\psi(\alpha)}_{n+1}} g^{-1}_{\xi^{\psi(\alpha)}_{n+2}})\cdots(g_{\xi^{\psi(\alpha)}_{n+m}} g^{-1}_{\xi^{\psi(\alpha)}_{n+m+1}})\right) \\
 & = & c(g_{\xi^{\psi(\alpha)}_n} g_{\xi^{\psi(\alpha)}_{n+m+1}}^{-1}) \\
 & = & d_{\psi(\alpha)}(\{\xi^{\psi(\alpha)}_n,\xi^{\psi(\alpha)}_{n+m+1}\}) \\
 & = & i_{\psi(\alpha)} \\
 & = & i.
\end{eqnarray*}
 \end{proof}

The following corollary gives a reasonably complete picture of which cardinals satisfy the Adjacent Hindman's Theorem, provided that their cofinality is larger than the number of colours.

\begin{corollary}
Let $\kappa, \lambda$ be cardinals. If either 
\begin{enumerate}
\item $\kappa < \omega$ and $\lambda  \geq \omega$; or
\item $\kappa \geq \omega$, $\lambda \geq (2^\kappa)^+$ and $\cf(\lambda)> \kappa$,
\end{enumerate}
then $\lambda \rightarrow (\lambda)^\afp_\kappa$.
\end{corollary}

\begin{proof}\hfill
\begin{enumerate}
\item If $\lambda = \aleph_0$ then we may use Hindman's Theorem, in its original countable version, to obtain, for any given colouring $c$, an infinite injective sequence $(g_n\big|n<\omega)$ such that all finite non-trivial products of distinct elements from that sequence have the same colour under $c$; {\em a fortiori} the desired result for finite {\em adjacent} products follows. If, on the other hand, $\lambda$ is uncountable, then the fact that $\lambda \rightarrow(\omega)^2_\kappa$ by Ramsey's Theorem allows us to use Lemma \ref{lem:ub}.
\item This is immediate from the Erd\H{o}s--Rado Theorem $\beth_1(\kappa)^+ \to (\kappa^+)^{2}_\kappa$ together with Lemma \ref{lem:ub}.
\end{enumerate}
\end{proof}

In view of the above results, it is natural to ask what happens when the cardinality of the group has cofinality 
no larger than the number of colours. We show that, while the Adjacent Hindman's Theorem fails if we require 
a monochromatic set the size of the group, we can ensure monochromatic sets of all smaller cardinalities.

\begin{theorem}
Let $\kappa, \lambda$ be cardinals, with $\lambda \geq (2^\kappa)^+$ and $\cf(\lambda)\leq \kappa$. If $G$ is a group of cardinality $\lambda$, then:
\begin{enumerate}
\item For each cardinal $\mu < \lambda$, $G \rightarrow(\mu)^\afp_\kappa$,
\item $G \nrightarrow(\lambda)^\afp_\kappa$.
\end{enumerate}
\end{theorem}

\begin{proof}\hfill
\begin{enumerate}
\item Let $c: G \longrightarrow \kappa$ be a colouring of $G$ with $\kappa$ colours. Upon fixing an independent sequence $(g_\alpha\big|\alpha<\lambda)$ in $G$ and a partition $C=\{C_\alpha\big|\alpha<\lambda\}$ of $\lambda$ into $\lambda$ pieces of cardinality $\lambda$, we proceed exactly as in the proof of Lemma \ref{lem:ub} to obtain, for each $\alpha < \lambda$, a colouring $d_\alpha:[C_\alpha]^2\longrightarrow \kappa$, a (strictly increasing) sequence $(\xi^\alpha_n \big | n < \omega)$ in $C_\alpha$, and a colour $i_\alpha < \kappa$ such that $[\{ \xi^\alpha_n \big | n < \omega\}]^2$ is $d_\alpha$-monochromatic with colour $i_\alpha$. This induces a partition of $\lambda$ into at most $\kappa$ cells $D_i = \{ \alpha < \kappa \big | i_\alpha = i\}$, 
for $i < \kappa$. 

Suppose that for each $i < \kappa$ the cardinality of $D_i$ is less than $\mu$. Then $|\lambda| \leq \kappa \mu$, which is impossible. Therefore, there exists an $i < \kappa$ such that $|D_i | \geq \mu$, i.e.
there are $\mu$-many $\alpha$'s less than $\kappa$ such that $i_\alpha = i$. Pick an injective function $\psi:\mu\longrightarrow D_i$ and set, for $\alpha<\mu$,
$x_{\omega \cdot \alpha + n} = g_{\xi^{\psi(\alpha)}_n} g^{-1}_{\xi^{\psi(\alpha)}_{n+1}}$ as in the proof of Lemma \ref{lem:ub}. The sequence $(x_\alpha \big | \alpha < \mu)$ is monochromatic for $c$ with colour $i$.
 
\item Since $\cf(\lambda) \leq \kappa$, there is a partition $(C_\alpha)_{\alpha < \cf(\lambda)}$ of $G$ such that $|C_\alpha|< \lambda$ for all $\alpha < \cf(\lambda)$. 
 Let $c:G \longrightarrow \kappa$ be the unique colouring satisfying $ c^{-1}[\{\alpha\}] = C_\alpha$ for all $\alpha < \cf(\lambda)$. No injective sequence of length $\lambda$ in $G$ is monochromatic for $c$. {\it A fortiori} there is no injective sequence of length $\lambda$ whose adjacent finite products are $c$-monochromatic.
\end{enumerate}
\end{proof}

\section{Failures of the uncountable Adjacent Hindman's Theorem}

So far we obtained complete information on whether $\aht(\lambda,\mu,\kappa)$ holds whenever $\lambda\geq(2^\kappa)^+$ (for infinite $\kappa$, regardless of the value of $\cf(\lambda)$). We now proceed to analyze the case $\lambda\leq 2^\kappa$. As we shall see, in this case one can find groups $G$ with $|G|=2^\kappa$ such that $G\nrightarrow(\mu)_\kappa^\afp$ even for very small values of $\mu$. In fact, we will show that $G\nrightarrow(2)_\kappa^\afp$ for all Abelian groups $G$ with $|G|=2^\kappa$; furthermore, we provide two particular examples of (non-Abelian) groups satisfying the same combinatorial properties (i.e., groups $G$ with $|G|=2^\kappa$ and $G\nrightarrow(2)_\kappa^\afp$). In other words, the lower bound with value $(2^\kappa)^+$ found in the previous section is optimal for the Adjacent Hindman's Theorem.

In light of the previous explanation, we state the theorem which is the main offshoot of the results from this section. Its proof can be taken to be any of the three particular results that follow.

\begin{theorem}\label{thm-main-optimal}
For any infinite cardinal $\kappa$, we have $2^\kappa\nrightarrow(2)_\kappa^\afp$.
\end{theorem}

The following three theorems provide different ways of establishing Theorem~\ref{thm-main-optimal}. The next theorem is, in fact, a generalization of~\cite[Theorem 10]{fernandez-lee}.

\begin{theorem}\label{thm:abelian-lower-bound}
If $G$ is an arbitrary Abelian group with $|G|=2^\kappa$, then $G\nrightarrow(2)_2^\afs$. More precisely, there exists a colouring $c:G\longrightarrow\kappa$ such that for all $x,y\in G$, $c(x)=c(y)$ implies $c(x)\neq c(x+y)$ (except in the trivial case where\footnote{Note that, according to our definitions, we only need to take care of injective sequences of elements of $G$ of length 2; that is, we would only need to prove that, if $x\neq y$ and $c(x)=c(y)$, then $c(x)\neq c(x+y)$. However, we can (and do) prove slightly more, since in our proof we really do not use the fact that $x\neq y$, but only that not both $x$ and $y$ equal $0$.} $x=y=0$).
\end{theorem}

\begin{proof}
If $G$ is an arbitrary Abelian group of size $2^\kappa$, then it is well-known (see e.g.~\cite[p. 123]{fernandez-nyjm}) that $G$ can be identified with a subgroup of $\bigoplus_{f\in 2^\kappa} G_f$, where each $G_f$ is equal to either $\mathbb Q$ or to the Pr\"ufer group $\mathbb Z[p^\infty]$, with $p$ a prime number. So assume without loss of generality that $G=\bigoplus_{f\in 2^\kappa} G_f$. Given an $x\in G\setminus\{0\}$, we consider its support, defined as follows: $\supp(x)=\{f\in 2^\kappa\big|x(f)\neq 0\}$. Given $f,g\in 2^\kappa$ we denote by $\Delta(f,g)=\min\{\alpha<\kappa\big|f(\alpha)\neq g(\alpha)\}$ (with the convention that $\Delta (f,f)=\kappa$). 
This allows us to define the following colouring, for $x\neq 0$:
\begin{equation*}
c(x)=
\begin{cases}
x(f)\ \ \ \ \ \ \ \ \ \ \ \ \ \ \ \text{ if }\supp(x)=\{f\}, \\
[\Delta(f_i,f_j)]_{i,j\leq n}\ \ \ \text{ if }\supp(x)=\{f_1,\ldots,f_n\} \\
\ \ \ \ \ \ \ \ \ \ \ \ \ \ \ \ \ \ \ \ \ \text{ with }n\geq 2\text{ and }f_1<\cdots<f_n.
\end{cases}
\end{equation*}
According to the definition, $c(x)$ is either an element of some $G_f$, or a doubly-indexed finite sequence---a matrix, if you will---of ordinals; in the latter case, note that $f_i<f_j$ is intended to refer to the usual lexicographic ordering on $2^\kappa$. We also let $c(0)=0$. Note that the range of $c$ is a set of size $\kappa$. Suppose by way of contradiction that there are non-zero elements $x,y\in G$ such that $c(x)=c(y)=c(x+y)$, and assume that we have taken such $x$ and $y$ with $|\supp(x)|$ as small as possible. If $\supp(x)=\{f\}$ then we must have $\supp(y)=\{f\}$ (otherwise $\supp(x)$ and $\supp(y)$ are disjoint and hence $c(x)\neq c(x+y)$) and so, letting $c(x)=x(f)=t=y(f)=c(y)$, we see that $c(x+y)=(x+y)(f)=2t\neq t=c(x)$, a contradiction (for otherwise, $2t=t$ would imply $t=0$, thus making $x=y=0$). Therefore, we may assume that $\supp(x)=\{f_1,\ldots,f_n\}$, $\supp(y)=\{g_1,\ldots,g_n\}$, and $\supp(x+y)=\{h_1,\ldots,h_n\}$, for some $n\geq 2$.

Let $\alpha$ be the least ordinal smaller than $\kappa$ that occurs as an entry of the matrix $[\Delta(f_i,f_j)]_{i,j\leq n}=[\Delta(g_i,g_j)]_{i,j\leq n}=[\Delta(h_i,h_j)]_{i,j\leq n}$; so there must be a single $f\in 2^\alpha$ with $f_i\upharpoonright\alpha=f$ for all $i$, and it must also be the case that $g_i\upharpoonright\alpha=f$ for all $i$ (else we would have $\supp(x)\cap\supp(y)=\varnothing$ and so $c(x)\neq c(x+y)$) and consequently also $h_i\upharpoonright\alpha=f$ for all $i$. Note that $f_1(\alpha)=0$ and $f_n(\alpha)=1$, so let $l\leq n$ be least such that $f_l(\alpha)=1$. Then $\Delta(f_1,f_i)=\alpha$ if and only if $i\geq l$; since $c(x)=c(y)=c(x+y)$ we may conclude that $l$ is also the least number such that $g_l(\alpha)=1$ and it is also the least number such that $h_l(\alpha)=1$. The conclusion is that, if we define $x'$ so that it agrees with $x$ on $\{f_1,\ldots,f_{l-1}\}$ and is zero everywhere else, and we analogously define $y'$ so that it agrees with $y$ on $\{g_1,\ldots,g_{l-1}\}$ and is zero everywhere else, then $x'+y'$ will agree with $x+y$ on $\{h_1,\ldots,h_{l-1}\}$ and be zero everywhere else; consequently, we will have $c(x')=c(y')=c(x'+y')$, contradicting the minimality of $|\supp(x)|$.

\end{proof}


Interestingly, we are also able to obtain non-Abelian groups witnessing the optimality of the lower bound from the previous section. Next comes the first such group.

\begin{definition}
Given a cardinal number $\kappa$, we let $S_\kappa$ be the {\bf symmetric group} on $\kappa$; in other words,
\begin{equation*}
S_\kappa=\{\sigma:\kappa\longrightarrow\kappa\big|\sigma\text{ is a bijection}\},
\end{equation*}
equipped with the composition of functions as the group operation.
\end{definition}

\begin{theorem}\label{permutation-group}
For any infinite cardinal $\kappa$, the statement $S_\kappa\nrightarrow(2)_\kappa^\afs$ holds; in other words, there exists a colouring $c:G\longrightarrow\kappa$ such that for $\sigma,\pi\in S_\kappa$, if $c(\sigma)=c(\pi)$ then $c(\sigma)\neq c(\sigma\pi)$ (except in the trivial case where $\sigma$ and $\pi$ are both the identity permutation).
\end{theorem}

\begin{proof}
Given a non-identity permutation $\sigma\in S_\kappa$ we let $\alpha_\sigma$ be the least ordinal moved by $\sigma$, and define $c(\sigma)=(\alpha_\sigma,\sigma(\alpha_\sigma))$; say also by convention that $c(\id_\kappa)=(0,0)$ (notice that the range of $c$ is a subset of $\kappa\times\kappa$, and so it has cardinality $\kappa$). Suppose we have two elements $\sigma,\pi\in S_\kappa$, not both $\id_\kappa$, with $c(\sigma)=c(\pi)=(\alpha,\beta)$. Then necessarily $\sigma\neq\id_\kappa\neq\pi$. Note that $\alpha_\sigma=\alpha_\pi=\alpha$, and every ordinal $\xi<\alpha$ is fixed by $\sigma$, $\pi$, and $\sigma\pi$. Now, there are two cases to consider:
\begin{description}
\item[Case 1] If $\sigma(\beta)=\alpha$, this means $\sigma(\pi(\alpha))=\alpha$ and so either $\sigma\pi$ is the identity permutation, or $\alpha_{\sigma\pi}>\alpha$. Either way, this implies that $c(\sigma\pi)\neq c(\sigma)$.
\item[Case 2] If $\sigma(\beta)=\gamma\neq\alpha$, then (trivially) also $\gamma\neq\beta$, and $\sigma(\pi(\alpha))=\gamma$. Hence $\alpha_{\sigma\pi}=\alpha$, and $c(\sigma\pi)=(\alpha,\gamma)\neq(\alpha,\beta)=c(\sigma)$.
\end{description}
The proof of Theorem~\ref{permutation-group} is finished.
\end{proof}

We now proceed to exhibit our second non-Abelian group of interest.

\begin{definition}
Let $\lambda$ be a cardinal, and let  $\{a_\alpha\big|\alpha<\lambda\}$ be a set of cardinality $\lambda$. The {\bf free group} on $\lambda$ generators, denoted by $F_\lambda$, consists of all reduced words on the alphabet $\{a_\alpha\big|\alpha<\lambda\}\cup\{a_\alpha^{-1}\big|\alpha<\lambda\}$ (a word is {\it reduced} if it contains no adjacent occurrences of the letters $a_\alpha$ and $a_\alpha^{-1}$ for any $\alpha$), with the group operation given by concatenation followed by reduction. The choice of $\{a_\alpha\big|\alpha<\lambda\}$ is immaterial.
\end{definition}

\begin{theorem}\label{free-group}
For any infinite cardinal $\kappa$, we have $F_{2^\kappa} \nrightarrow (2)^\afp_\kappa$.
\end{theorem}

\begin{proof}
Let $\{ a_f \big | f : \kappa \longrightarrow 2\}$ be the set of free generators of $F_{2^\kappa}$. Note that each non-identity element $w\in F_{2^\kappa}$ can be uniquely expressed, for some $k\in\mathbb N$, as $w=a_{f_1}^{n_1} \cdots a_{f_k}^{n_k}$, 
where $f_i\neq f_{i+1}$ for each $i\in\{1,\ldots,k-1\}$, and each $n_i$ is in $\mathbb{Z}\setminus \{0\}$. We will say that $k$ is the {\it reduced length} of the word $w$, written $k=\ell(w)$.

With the notion of reduced length at our disposal, we define the colouring $c$ as follows:
\begin{equation*}
c(w)=
\begin{cases}
0\ \ \ \ \ \ \ \ \ \ \ \ \ \ \ \ \ \ \ \ \ \ \ \ \ \ \ \ \text{ if } w= e\text{ (the empty word)}, \\
(\ell(w),\Delta(f_k,f_{k+1}))\ \ \ \ \ \text{ if }\ell(w)=2k\text{ and }w =a_{f_1}^{n_1}\cdots a_{f_{2k}}^{n_{2k}}, \\
(\ell(w),n_{k+1})\ \ \ \ \ \ \ \ \ \ \ \ \ \ \text{ if } \ell(w)=2k+1\text{ and }w =a_{f_1}^{n_1}\cdots a_{f_{2k+1}}^{n_{2k+1}},
\end{cases}
\end{equation*}
where, as in the proof of Theorem~\ref{thm:abelian-lower-bound}, $\Delta(f,g)$ denotes the least $\alpha<\kappa$ such that $f(\alpha)\neq g(\alpha)$ for two distinct $f,g\in 2^\kappa$.
Note that $c:F_{2^\kappa}\longrightarrow\{0\}\cup(\mathbb N\times(\mathbb Z\setminus\{0\}))\cup(\mathbb N\times\kappa)$ and so $c$ is a colouring with at most $\kappa$ colours. Crucially, any two words with the same colour must have the same reduced length. So suppose that $w,v\in F_{2^\kappa}$ are not both the empty word, and $c(w)=c(v)$. Then, indeed, neither $w$ nor $v$ is the empty word, so $\ell(w)=\ell(v)\geq 1$. Again our argument splits as follows:

\begin{description}
\item[Case 1] Suppose that $\ell(w)=\ell(v)$ is even, say $2k$, so we may write $c(w)=c(v)=(2k,\alpha)$ for some $\alpha<\kappa$. Write $w=a_{f_1}^{n_1}\cdots a_{f_{2k}}^{n_{2k}}$ and $v=a_{g_1}^{m_1}\cdots a_{g_{2k}}^{m_{2k}}$ with $f_i\neq f_{i+1}$ and $g_i\neq g_{i+1}$ for each $i$, and $n_1,\ldots,n_{2k}$, $m_1,\ldots,m_{2k}\in\mathbb Z\setminus\{0\}$. If $\ell(wv)\neq\ell(w)$ then $c(wv)\neq c(w)$, so we may assume $\ell(wv)=2k$ and note that this can only happen if $f_k\neq g_{k+1}$ and the subwords (of $w$ and $v$ respectively) $a_{f_{k+1}}^{n_{k+1}}\cdots a_{f_{2k}}^{n_{2k}}$ and $a_{g_1}^{m_1}\cdots a_{g_k}^{m_k}$ are inverses of each other. In particular, $a_{f_{k+1}}^{n_{k+1}}=(a_{g_k}^{m_k})^{-1}$, so $f_{k+1}=g_k$. Furthermore, $\alpha=\Delta(f_k,f_{k+1})=\Delta(g_k,g_{k+1})$, so $\alpha=\Delta(f_{k+1},g_{k+1})$. Hence the two (distinct) functions $f_k$ and $g_{k+1}$ agree with $f_{k+1}$ up to any $\xi<\alpha$ but disagree with it at $\alpha$. This implies that $\Delta(f_k,g_{k+1})>\alpha$. We thus have
\begin{equation*}
c(wv)=c(a_{f_1}^{n_1}\cdots a_{f_k}^{n_k} a_{g_{k+1}}^{m_{k+1}}\cdots a_{g_{2k}}^{m_{2k}})=(2k,\Delta(f_k,g_{k+1}))\neq(2k,\alpha)=c(w),
\end{equation*}
and we are done.

\item[Case 2] Now suppose that $\ell(w)=\ell(v)$ is odd, say $c(w)=c(v)=(2k+1,n)$ for some $n\in\mathbb Z\setminus\{0\}$ and some $k\in\omega$. Write\newline $w=a_{f_1}^{n_1}\cdots a_{f_k}^{n_k} a_{f_{k+1}}^n a_{f_{k+2}}^{n_{k+2}}\cdots a_{f_{2k+1}}^{n_{2k+1}}$ and $v=a_{g_1}^{m_1}\cdots a_{g_k}^{m_k} a_{g_{k+1}}^n a_{g_{k+2}}^{m_{k+2}}\cdots a_{g_{2k+1}}^{m_{2k+1}}$ with $f_i\neq f_{i+1}$ and $g_i\neq g_{i+1}$ for each $i$, and $n_1,\ldots,n_{2k+1},m_1,\ldots,m_{2k+1}\in\mathbb Z\setminus\{0\}$. As in the previous case, we may assume without loss of generality that $\ell(wv)=2k+1$ (otherwise $c(wv)\neq c(w)$ and we are done). Note that, in order to have $\ell(wv)=2k+1$, we must have $f_{k+1}=g_{k+1}$ and the subword $a_{f_{k+2}}^{n_{k+2}}\cdots a_{f_{2k+1}}^{n_{2k+1}}$ (of $w$) must be the inverse of the subword $a_{g_1}^{m_1}\cdots a_{g_k}^{m_k}$ (of $v$). Hence,
\begin{eqnarray*}
c(wv)  & = & c(a_{f_1}^{n_1}\cdots a_{f_k}^{n_k}a_{f_{k+1}}^{2n} a_{g_{k+2}}^{m_{k+2}}\cdots a_{g_{2k+1}}^{m_{2k+1}}) \\
 & = & (2k+1,2n)\neq(2k+1,n)=c(w),
\end{eqnarray*}
which finishes the proof of Theorem~\ref{free-group}.
\end{description}
\end{proof}

\section*{Acknowledgements}

The second author was partially supported, at the beginning of the present work, by a DGAPA postdoctoral fellowship from Universidad Nacional Aut\'onoma de M\'exico; and by grant SIP-20221862 from the Instituto Polit\'ecnico Nacional towards the end of it.

\end{document}